\documentclass[11pt]{article}
\usepackage[T1]{fontenc}
\usepackage{latexsym,amssymb,amsmath,amsfonts,amsthm}
\newcommand{\R}{{\mathbb R}}
\newcommand{\Z}{{\mathbb Z}}
\newcommand{\N}{{\mathbb N}}
\newcommand{\C}{{\mathbb C}}
\newcommand{\Rc}{{\mathcal R}}
\newcommand{\ds}{\displaystyle}
\newcommand{\no}{\nonumber}
\newcommand{\be}{\begin{eqnarray}}
\newcommand{\ben}{\begin{eqnarray*}}
\newcommand{\en}{\end{eqnarray}}
\newcommand{\enn}{\end{eqnarray*}}
\newcommand{\ba}{\backslash}
\newcommand{\bt}{\beta}

\newcommand{\curl}{{\rm curl\,}}
\newcommand{\dive}{{\rm div\,}}
\newcommand{\I}{{\rm Im\,}}
\newcommand{\Rt}{{\rm Re\,}}

\newcommand{\g}{\gamma}
\newcommand{\f}{\frac}
\newcommand{\G}{\Gamma}

\newcommand{\Om}{\Omega}
\newcommand{\om}{\omega}

\newcommand{\al}{\alpha}

\newtheorem{theorem}{Theorem}[section]
\newtheorem{lemma}[theorem]{Lemma}

\begin{document}
\renewcommand{\theequation}{\arabic{section}.\arabic{equation}}
\begin{titlepage}
\title{\bf An inverse electromagnetic scattering problem for a
bi-periodic inhomogeneous layer on a perfectly conducting plate}
\author{Guanghui Hu\\
Weierstrass Institute for Applied Analysis and Stochastics\\
Mohrenstr. 39, 10117 Berlin, Germany\\ \\
Jiaqing Yang and Bo Zhang\\
LSEC and Institute of Applied Mathematics\\
Academy of Mathematics and Systems Science\\
Chinese Academy of Sciences\\
Beijing 100190, China\\
({\sf hu@wias-berlin.de} (GH),\ {\sf jiaqingyang@amss.ac.cn}(JY),\ {\sf b.zhang@amt.ac.cn} (BZ))}

\date{}
\end{titlepage}
\maketitle

%\vspace{.2in}
%\noindent Proposed running head:
%{\bf Inverse electromagnetic scattering for bi-periodic inhomogeneous media}\\
%\vspace{.1in}
%
%\noindent
%Manuscript Correspondence Address:\\
%Professor Bo Zhang\\
%Institute of Applied Mathematics\\
%Academy of Mathematics and Systems Science\\
%Chinese Academy of Sciences\\
%Beijing 100190, China\\
%Email: b.zhang@amt.ac.cn\\
%Telephone: +86 10 6265 1358\\
%Fax:       +86 10 6254 1689
%
%\newpage
%\begin{center}
%\title{\Large\bf An inverse electromagnetic scattering problem for a
%periodic inhomogeneous layer on a perfectly conducting plate}
%\end{center}
%
%\vspace{.2in}

\begin{abstract}
This paper is concerned with uniqueness for reconstructing a periodic
inhomogeneous medium covered on a perfectly conducting plate. We deal with the problem
in the frame of time-harmonic Maxwell systems without $TE$ or $TM$ polarization.
An orthogonal relation for two refractive index is obtained, and then inspired by Kirsch's idea,
the refractive index can be identified by utilizing the eigenvalues and
eigenfunctions to a quasi-periodic Sturm-Liouville eigenvalue problem.

\vspace{.2in}
\end{abstract}

\section{Introduction}
\setcounter{equation}{0}

Scattering theory in periodic structures has many applications in micro-optics,
radar imaging and nondestructive testing. We refer to \cite{PR} for historical
remarks and details of these applications.
Consider a time-harmonic electromagnetic plane wave incident on a bi-periodic
layer sitting on a perfectly conducting plate in $\R^3$.
We assume that the medium inside the layer consists of some inhomogeneous
isotropic conducting or dielectric material, whereas the medium above the layer
consists of some homogeneous dielectric material. Suppose the magnetic
permeability is a fixed positive constant throughout the whole space.
The material properties of the media are then characterized completely
by an index of refraction in the layer and a positive constant above the layer.
The direct scattering problem is, given the incident field and the bi-periodic
refractive index, to study the electromagnetic distributions, whereas the
inverse scattering problem is to determine the refractive index from the
knowledge of the incident waves and their corresponding measured scattered fields.

Adopting the Cartesian axis $ox_1x_2x_3$ with the $x_3$-axis vertically upwards,
perpendicular to the plate. If the refractive index is invariant in the $x_2$
direction, the direct and inverse problems as indicated above can be dealt
with in the TE polarization case where the electric field $E(x)$ is transversal to
the $(x_1, x_3)$-plane by assuming $E=(0,u(x_1,x_3),0)$, or in the TM polarization
case where the magnetic field $H(x)$ is transversal to the $(x_1, x_3)$-plane
by assuming $H=(0,u(x_1,x_3),0)$. In the case of TE polarization,
Kirsch \cite{Kirsch95} has studied the direct scattering problem via the variational
method, and for the inverse problem, instead of constructing the complex
geometrical optical solutions as in the Calder\'{o}ns problem
(see \cite{Nachman,Uhlmann}), he considered a class of eigenfunctions
to a special kind of quasi-periodic Sturm-Liouville eigenvalue problem.
Relying on the asymptotic behavior of those eigenvalues,
the uniqueness result for the inverse problem can be proved once the orthogonal
relation for two different refractive indexes has obtained.
See also \cite{S98,S99} for the direct and inverse acoustic
scattering by periodic, inhomogeneous, penetrable medium in the whole $\R^2$.
Other uniqueness results for reconstructing the profile of a  bi-periodic
perfectly conducting grating can be seen in \cite{AH,BZ,BZZ09}.

In this paper, we are mainly concerned with the uniqueness issue for
reconstructing the refractive index in the framework of time-harmonic Maxwell equations
without TE or TM polarization. The uniqueness result for the inverse problem
in this paper is most closely related in term of result and method of argument
to Kirsch on the determination of the refractive index in the TE polarization.
Inspired by \cite{SU} and \cite{Iaskov}, we obtain an orthogonality relation for
two different refractive indexes by using a D-to-N map on an artificial boundary on
which the tangential electric fields are identical for an integral type of incident
electric field. It should be remarked that the method for constructing geometry optical
solutions in \cite{Nachman,Iaskov,SU} for non-periodic inverse conductivity problems
does not work since the solutions are required to be quasi-periodic in the periodic case.
To reconstruct the refractive index, we follow Kirsch's idea \cite{Kirsch95}
(see also \cite{S99}) by considering a kind of Sturm-Liouville eigenvalue problems.
We shall prove the uniqueness result when the index depends only on one direction
($x_1$ or $x_2$). However, we expect the result to hold in a more general case
by constructing special solutions with suitable asymptotic behaviors for the Maxwell
equations.
% which needs more techniques.
%As far as the authors' knowledge, this is the first paper dealing with the
%uniqueness for reconstructing the periodic inhomogeneous refractive index
%in the frame of Maxwell system.

%Since we focus on the inverse problem, we only give a sketch of proof on the existence
%and uniqueness of solutions to the direct problem which have received an
%increasing studies in recent years
Scattering by bi-periodic structures have been studied by many authors using both integral
equation methods and variational methods (see, e.g. \cite{AT}, \cite{B}, \cite{D},
\cite{DF}, \cite{EG98}, \cite{Kirsch1}, \cite{NS}  and \cite{Schmidt}).
It is known that, for all but possibly a discrete set of frequencies,
the direct scattering problem has a unique weak solution in the case of
bi-periodic inhomogeneous medium in the whole $\R^3$, of which an absorbing medium
always leads to a uniqueness result for any frequency. When the refractive
index is non-absorbing, uniqueness can be guaranteed in the TE mode if the refractive index
satisfies an increasing criterion in the $x_3$-direction (\cite{S98,BB94}).
See also \cite{Zhang1} and \cite{Zhang2} for the uniqueness results of
more general rough surface scattering by an inhomogeneous medium in a half space
in the TE or TM mode. In this paper, we assume that the medium inside the layer is
absorbing so that the uniqueness result for the direct problem holds, implying that
the D-to-N map $T$ (at the end of Section \ref{sec3}), which depends on the refractive index,
is well-defined.

The rest of the paper is organized as follows. In the next section we set up
the precise mathematical framework and introduce some quasi-periodic function spaces needed.
In Section \ref{sec3}, we consider a quasi-periodic boundary value problem (QPBVP)
in a periodic cell via the variational approach which is used for the study of the inverse problem.
Uniqueness and existence of solutions to the QPBVP are justified by the classic Hodge
decomposition and the Fredholm alternative. This leads to the definition of a D-to-N map
on an artificial boundary which is continuous and depends on the refractive index.
In Section \ref{sec4}, based on the property of the transparent boundary
condition defined on the artificial boundary, we give a solvability result
of the direct scattering problem. In Section \ref{inverse problems}, we establish a
uniqueness result for the inverse scattering problem.

\section{Time-harmonic Maxwell equations and quasi-periodic function spaces}
\setcounter{equation}{0}

\subsection{Time-harmonic Maxwell equations}

Let $\mathbb{R}^3_+=\{(x_1,x_2,x_3)\in\R^3\,|\,x_3>0\}$ and assume that $\mathbb{R}^3_+$
is filled with an inhomogeneous, isotropic, conducting or dielectric medium of electric
permittivity $\epsilon>0$, magnetic permeability $\mu>0$ and electric conductivity $\sigma\geq 0$.
Suppose the medium is non-magnetic, that is, the magnetic permeability $\mu$ is a fixed constant in
$\mathbb{R}^3_+$ and the field is source free. Then the electromagnetic wave propagation is governed
by the time-harmonic Maxwell equations (with the time variation of the form $e^{-i\om t},$ $\om>0$)
\be\label{1.1}
 \curl E-i\omega\mu H=0,\qquad\;\;\curl H+i\omega(\epsilon+i\f{\sigma}{\omega})E=0,
\en
where $E$ and $H$ are the electric field and magnetic field, respectively. Suppose
the inhomogeneous medium is $2\pi$-periodic with respect to $x_1$ and $x_2$ directions,
that is, for all $n_1,n_2\in\mathbb{Z}^2$,
\ben
\epsilon(x_1+2\pi n_1,x_2+2\pi n_2,x_3)=\epsilon(x_1,x_2,x_3),\qquad
\sigma(x_1+2\pi n_1,x_2+2\pi n_2,x_3)=\sigma(x_1,x_2,x_3).
\enn
Further, assume that
%there is some $\delta>0$ such that
$\epsilon(x)=\epsilon_0,$  $\sigma=0$ for $x_3>b$ (which means that the medium above
the layer is lossless) and that the inhomogeneous medium has a perfectly conducting
boundary $\G_0:=\{x_3=0\}$. Consider a time-harmonic plane wave
\ben
E^i=pe^{ikx\cdot d},\qquad  H^i=qe^{ikx\cdot d},
\enn
incident on the periodic inhomogeneous layer from the top region $\Omega:=\{x\in\R^3\,|\,x_3>b\}$,
where $d=(\al_1,\al_2,-\bt)=(\cos\theta_1\cos\theta_2,\cos\theta_1\sin\theta_2,-\sin\theta_1)$
is the incident wave vector whose direction is specified by $\theta_1$ and $\theta_2$ with
$0<\theta_1<\pi,\,0<\theta_2\leq 2\pi$ and the vectors $p$ and $q$ are polarization directions
satisfying that $p=\sqrt{{\mu}/{\varepsilon}}(q\times d)$ and $q\bot d.$
The problem of scattering of time-harmonic electromagnetic waves in this model leads to the
following problem:
\be\label{equation1}
 \curl\curl E-k^2E=0  &&\mbox{in}\quad x_3>b, \\ \label{equation2}
 \curl\curl E-k^2qE=0  &&\mbox{in}\quad\Om_b,\\ \label{equation3}
 \nu\times E=0       &&\mbox{on}\quad\G_0,\\ \label{equation4}
 E=E^i+E^s  &&\mbox{in}\quad\mathbb{R}^3_+,
\en
where $k=\sqrt{\epsilon_0\mu}\omega$ is the wave number,
$q(x)=\frac{1}{\epsilon_0}(\epsilon(x)+i\frac{\sigma(x)}{\omega})$ is the
refractive index and $\nu$ is the unit normal at the boundary.

Set $\al=(\al_1,\al_2,0)\in\R^3$ and $n=(n_1,n_2)\in\Z^2$. The periodicity of the medium
motivates us to look for $\al$-quasi-periodic solutions in the sense that
$E(x_1,x_2,x_3)e^{-i\al\cdot x}$ is $2\pi$ periodic with respect to $x_1$ and $x_2$,
respectively. Since the domain is unbounded in the $x_3$-direction, a radiation condition must be
imposed. It is required physically that the diffracted fields remain bounded as $x_3$
tends to $+\infty$, which leads to the so-called outgoing wave condition in the form of
\be\label{RE}
E^s(x)=\ds\sum_{n\in\Z^2}E_ne^{i(\al_n\cdot x+\bt_nx_3)},\qquad x_3>b,
\en
where $\al_n=(\al_1+n_1,\al_2+n_2,0)\in\R^3$, $E_n=(E_n^{(1)},E_n^{(2)},E_n^{(3)})\in\C^3$
are constant vectors and
$$
\beta_n=\left\{\begin{array}{lll}
             (k^2-|\al_n|^2)^{\frac{1}{2}}\qquad\rm{if}\ |\al_n|< k,\\
              i(|\al_n|^2-k^2)^{\frac{1}{2}}\qquad\rm{if}\ |\al_n|> k,\\
             \end{array}\right.
$$
with $i^2=-1$. Furthermore, we assume that $\beta_n\neq 0$ for all $n\in\Z^2$.
The series expansion in (\ref{RE}) is considered as the Rayleigh series of the scattered
field and the condition is called the Rayleigh expansion radiation condition.
The coefficients $E_n$ in (\ref{RE}) are also called the Rayleigh sequence.
From the fact that $\rm{div}\,E^s(x)=0$ it is clear that
\ben
   \al_n\cdot E_n+ \beta_nE_n^{(3)}=0.
\enn

The direct problem (DP) is to compute the total field $E$ in $\mathbb{R}^3_+$, given the
incident wave $E^i$, the refractive index $q(x)$ and the boundary condition on $\G_0$.
Since only a finite number of terms in (\ref{RE}) are upward propagating plane waves
and the rest is evanescent modes that decay exponentially with distance away
from the periodic medium, we use the near field data rather than the far field data to
reconstruct the refractive index $q(x)$.
Thus, our inverse problem (IP) is to determine the periodic medium $q(x)$ from
a knowledge of the incident wave $E^i$ and the total tangential electric
field $\nu\times E$ on a plane $\G_a=\{x\in\mathbb{R}^3\,|\,x_3=a\}(a>b)$ above the layer.

%{\em Scattering of electromagnetic waves by a smooth doubly periodic structure
%has been studied by many authors using both integral and variational methods.
%See, e.g. \cite{AT,AH,B,D,DF,NS} for results on existence,
%uniqueness, and numerical approximations of solutions to the direct problems.
%The inverse problem in a homogeneous medium with a perfect conductor under the
%doubly periodic structure has been considered in \cite{AH,BZ}. With a lossy
%homogenous medium (i.e., $\I(k)>0$) above the conductor, Ammari \cite{AH} proved
%a global uniqueness result for the inverse problem with one incident plane wave.
%For the case of lossless medium (i.e., $\I(k)=0$) above the conductor, a local
%uniqueness result was obtained by Bao and Zhou in \cite{BZ} for the inverse problem
%with one incident plane wave by establishing a lower bound
%of the first eigenvalue of the $\curl\curl$ operator with the boundary condition
%(\ref{equation2}) in a bounded, smooth convex domain in $\R^3.$ The stability of
%the inverse problem was also studied in \cite{BZ}.
%For inverse scattering problems by bounded obstacles the reader is referred to \cite{CK}.}

\subsection{Quasi-periodic function spaces}

In this section we introduce some function spaces needed for the scattering problem
(\ref{equation1})-(\ref{equation4}).
These spaces will play a crucial role not only in the study of the direct problem
but also in the inverse problem.
In \cite{B,D,Schmidt}, the authors always seek the $H^1$-variational solution for
the magnetic field $H$, based on the facts that the magnetic permeability $\mu>0$
is a constant and that any vector field $H\in L^2{(D)}^3$ satisfying
that $\triangledown\times H\in L^2{(D)}^3$ and $\triangledown\cdot H\in L^2{(D)}^3$
belongs to $H^1_{loc}(D)^3$ for any bounded domain $D\subset\R^3$.
In this paper, based on the classic Hodge decomposition, we are interested in
weak solutions in $H(\curl)$ of the problem (\ref{equation1})-(\ref{equation4}), that is,
both $E$ and $\triangledown\times E$ belong to $L^2_{loc}(\R^3_+)^3$.
This allows us to solve the scattering problem in a general case
when $\mu$ is a periodic variable function other than a constant.

The scattering problem can be reduced to a single periodic cell. To this end, we reformulate
the following notations.
\ben
\G_b=\{x_3=b\,|\,0<x_1,x_2<2\pi\},\
\Om_b=\{x\in\R^3_+\,|\, x_3<b,\;0<x_1,x_2<2\pi\}.
\enn
We also need the following scalar quasi-periodic Sobolev space:
\ben
H^1(\Om_b)=\{u(x)=\sum_{n\in\Z^2}u_n(x_3)\exp(i\al_n\cdot x)\,|\,
u\in L^2(\Om_b),\nabla u\in (L^2(\Om_b))^3,u_n\in\C\}.
\enn
Denote by $H^{\frac12}(\G_b)$ the trace space of $H^1(\Om_b)$ on $\G_b$ with the norm
\ben
||f||^2_{H^{\frac12}(\G_b)}=\sum_{n\in\Z^2}|f_n|^2(1+|\al_n|^2)^{\frac12},\qquad f\in H^{\frac12}(\G_b),
\enn
where $f_n=(f,\exp(i\al_n\cdot x))_{L^2(\G_b)}$
and write $H^{-\frac12}(\G_b)=(H^{\frac{1}{2}}(\G_b))^\prime$, the dual space to $H^{\frac12}(\G_b)$.

We now introduce some vector spaces. Let
\ben
H(\curl,\Om_b)&=&\{E(x)=\sum_{n\in\Z^2}E_n(x_3)\exp(i\al_n\cdot x)\,|\,E_n\in\C^3,\\
              &&\quad\,E\in(L^2(\Om_b))^3,\,\curl E\in(L^2(\Om_b))^3\}
\enn
with the norm
\ben
||E||^2_{H(\curl,\Om_b)}=||E||^2_{L^2(\Om_b)}+||\curl E||^2_{L^2(\Om_b)}.
\enn
For $x^\prime=(x_1,x_2,b)\in\G_b$, $s\in\R$ define
\ben
H_t^s(\G_b)&=&\{E(x^\prime)=\sum_{n\in\Z^2}E_n\exp(i\al_n\cdot x^\prime)\,|\,
              E_n\in\C^3,\,e_3\cdot E=0,\\
           &&\qquad \|E\|^2_{H^s(\G_b)}=\sum_{n\in\Z^2}(1+|\al_n|^2)^s|E_n|^2<+\infty\}\\
H_t^s(\dive,\G_b)&=&\{E(x^\prime)=\sum_{n\in\Z^2}E_n\exp(i\al_n\cdot x^\prime)\,|\,
                     E_n\in\C^3,\,e_3\cdot E=0,\\
                 &&\quad ||E||^2_{H^s(\dive,\G_b)}=\sum_{n\in\Z^2}(1+|\al_n|^2)^s
                        (|E_n|^2+|E_n\cdot\al_n|^2)<+\infty\}\\
H_t^s(\curl,\G_b)&=&\{E(x^\prime)=\sum_{n\in\Z^2}E_n\exp(i\al_n\cdot x^\prime)\,|\,
                     E_n\in\C^3,\,e_3\cdot E=0,\\
                 &&\quad ||E||^2_{H^s(\curl,\G_b)}=\sum_{n\in\Z^2}(1+|\al_n|^2)^s
                     (|E_n|^2+|E_n\times\al_n|^2)<+\infty\}
\enn
and write $L_t^2(\G_b)=H_t^0(\G_b).$ Recall that
\ben
H_t^{-{1}/{2}}(\dive,\G_b)=\{e_3\times E|_{\G_b}\,|\,E\in H(\curl,\Om_b)\}
\enn
and that the trace mapping from $H(\curl,\Om_b)$ to $H_t^{-{1}/{2}}(\dive,\G_b)$ is continuous and
surjective (see \cite{BCS} and the references there).

It is well-known (see \cite{NS}) that the free space $\al$-quasi-periodic Green function
for the Helmholtz equation $(\Delta+k^2)u=0$ in $\R^3$ is given by
\be\label{Green}
G(x,y)=\frac{1}{8\pi^2}\sum_{n\in\Z^2}\frac{1}{i\beta_n}\exp(i\al_n\cdot(x-y)+i\beta_n|x_3-y_3|).
\en

We assume throughout this paper that $q$ satisfies the following conditions:
\begin{description}
\item ({\bf A1}) $q\in C^1(\overline{\Om_b})$ and $q(x)=1$ when $x_3>b$;
\item ({\bf A2}) $\I[q(x)]\geq0$ for all $x\in\overline{\Om_b}$ and
                 $\I[q(x_0)]>0$ for some $x_0\in\overline{\Om_b}$;
\item ({\bf A3}) $\Rt[q(x)]\ge\g$ for all $x\in\overline{\Om_b}$ for some positive constant $\gamma$.
\end{description}

\section{A quasi-periodic boundary value problem}\label{sec3}
\setcounter{equation}{0}

Before studying the original problem (\ref{equation1})-(\ref{RE}),
we first consider the following quasi-periodic boundary value problem in $\Om_b$:
\be\label{BVP}
  \curl\curl E-k^2q(x)E&=&0\qquad{\rm in}\;\Om_b,\\ \label{BV1}
  \nu\times E&=&0 \qquad {\rm on}\;\G_0,\\\label{BV2}
  \nu\times E&=&f \qquad {\rm on}\;\G_b,
\en
where $f\in H^{-1/2}_{\dive}(\G_b)$.

\begin{lemma}\label{wellpose1}
If the conditions $(A1)-(A3)$ are satisfied, then there exists a unique solution
$E\in H(\curl,\Om_b)$ to the problem $(\ref{BVP})-(\ref{BV2})$ such that
\ben\label{estimate}
||E||_{H(\curl,\Om_b)}\leq C||f||_{H^{-1/2}_{\dive}(\G_b)},
\enn
where $C$ is a positive constant independent of $f$.
\end{lemma}

\begin{proof}
We first prove the uniqueness part.
Let $f=0.$ Multiplying both sides of $(\ref{BVP})$ by $\overline{E}$
it follows from Green's vector formula, the quasi-periodic property of $E$ and
the boundary conditions (\ref{BV1}) and (\ref{BV2}) that
\be
\int_{\Om_b}[|\curl E|^2-k^2q|E|^2]dx=0.
\en
Take the imaginary part of the above equation and use the assumption on $q(x)$
to find that
\ben
\int_{B_{\epsilon}(x_0)}|E(x)|^2dx=0,
\enn
where $B_{\epsilon}(x_0)\subset\Om_b$ is a small ball centered at $x_0$ with radius $\epsilon$.
Thus $E(x)\equiv 0$ in $B_{\epsilon}(x_0)$.
By \cite[Theorem 6]{C96} we have $E\in (H^1(\Om_b))^3$. Thus, by the unique continuation
principle (see \cite[Theorem 2.3]{Ok02})
%By the unique continuation principle (see \cite[Lemma 8]{CK})
%%which implies that $\nu\times E=0,\qquad\nu\times\curl E=0$ on a regular surface $\G'$
%%contained in $B_{\epsilon}(x_0)$. By Helmgren's uniqueness theorem,
we have $E\equiv0$ in $\Om_b$.

We are now in a position to prove the existence of solutions.
For any $V\in H(\curl,\Om_b)$ such that $\nu\times E=0$ on $\G_0\cup\G_b$,
multiplying both sides of $(\ref{BVP})$ by $\overline{V}$ yields
\be\label{variational1}
\int_{\Om_b}[\curl E\cdot\curl \overline{V}-k^2q E\cdot\overline{V}]dx=0.
\en
There exists at least one element $W\in H(\curl,\Om_b)$ satisfying that
$\nu\times W=0$ on $\G_0$ and $\nu\times W=f$ on $\G_b$.
%Both sides of (\ref{variational1}) subtract the following quantity
%\ben
%\int_{\Om_b}\curl W\cdot\curl \overline{V}-k^2q W\cdot \overline{V}dx,
%\enn
Then the equation (\ref{variational1}) can be rewritten as
\ben
\int_{\Om_b}[\curl(E-W)\cdot\curl\overline{V}-k^2q(E-W)\cdot\overline{V}]dx
=-\int_{\Om_b}[\curl W\cdot\curl\overline{V}-k^2q W\cdot\overline{V}]dx.
\enn
Let $X:=\{U\in H(\curl,\Om_b),\;\nu\times U=0\;\mbox{on}\;\G_0\cup\G_b\}$.
Then $U:=E-W\in X$. Thus the problem (\ref{BVP})-(\ref{BV2}) is equivalent to
the following variational problem: Find $U\in X$ such that for any $V\in X$,
\be\label{variationl}
\int_{\Om_b}[\curl U\cdot\curl\overline{V}-k^2q U\cdot\overline{V}]dx=F_W(V),
\en
where $\ds F_W(V)=-\int_{\Om_b}[\curl W\cdot\curl\overline{V}-k^2q W\cdot\overline{V}]dx$.
The proof is broken down into the following steps.

{\bf Step 1.} To establish the Hodge decomposition:
\be\label{H-decom}
X=X_0\oplus\nabla S,
\en
where $S=\{p\in H^1(\Om_b),\;p=0\;{\rm on}\;\G_0\cup\G_b\}$ and
$\ds X_0=\{\xi\in X\,|\,\int_{\Om_b}q(x)\xi\cdot\nabla\overline{p}=0,\;\forall p\in S\}$.

For $U,V\in X$ define
$$
a(U,V)=\int_{\Om_b}[\curl U\cdot\curl\overline{V}-k^2q U\cdot\overline{V}]dx.
$$
It follows from the assumptions (A1)-(A3) on $q(x)$ that
\ben
|a(\nabla p,\nabla p)|\ge k^2\int_{\Om_b}\Rt[q(x)]|\nabla p|^2dx
\ge k^2\gamma||\nabla p||_{L^2(\Om_b)}^2=k^2\gamma||\nabla p||_{H(\curl,\Om_b)}^2.
\enn
Thus, for every $E\in X$ there exits a unique $p\in S$ such that
$a(\nabla p,\nabla q)=a(E,\nabla q)$ for all $q\in S$.
Let $\xi:=E-\nabla p$. Then it is easy to show that $\xi\in X_0$ and $X_0\cap S=\emptyset$,
which implies the Hodge decomposition (\ref{H-decom}).

{\bf Step 2.} To prove the existence of a unique solution $U\in X$
to the problem ((\ref{variationl}).
%such that $a(U,V)=F_W(V)$ for all $V\in X$.

By (\ref{H-decom}) we may assume that $U=\xi+\nabla p,\;V=\eta+\nabla q$
with $\xi,\eta\in X_0$ and $p,q\in S.$ Then the problem (\ref{variationl})
becomes the following one: Find $\xi\in X_0$ and $p\in S$ such that
\ben\label{variational2}
a(\nabla p,\nabla q)+a(\xi,\eta)=F_W(\nabla q)+F_W(\eta).
\enn
Since $a(\cdot,\cdot)$ is coercive on $\nabla S$, there exists a unique $p\in S$ such that
\ben
a(\nabla p,\nabla q)=F_W(\nabla q)\qquad \forall\,q\in S
\enn
with the estimate $\ds||\nabla p||_{H(curl,\Om_b)}\leq C||W||_{H(curl,\Om_b)}$.
It remains to find $\xi\in X_0$ such that $a(\xi,\eta)=F_W(\eta)$ for all $\eta\in X_0$.
The bilinear form $a(\cdot,\cdot)$ can be decomposed into the sum of the following two forms:
\ben
a_1(\xi,\eta)&=&\int_{\Om_b}\curl\xi\cdot\curl \overline{\eta}+\xi\cdot\overline{\eta}dx,\\
a_2(\xi,\eta)&=&-k^2\int_{\Om_b}(1+q)\xi\cdot\overline{\eta}dx.
\enn
Obviously, $a_1(\cdot,\cdot)$ is coercive on $X_0$, and
it follows from \cite[Lemma 3.2]{Habib} that $X_0$ is compactly imbedded into
$(L^2(\Om_b))^3$. Thus, by the standard Fredholm alternative theory there exists
a unique $\xi\in X_0$ satisfying that $a(\xi,\eta)=F_W(\eta)$ for all $\eta\in X_0$.
Furthermore, $\ds||\xi||_{H(curl,\Om_b)}\leq C||W||_{H(curl,\Om_b)}$.

{\bf Step 3.} To establish the estimate (\ref{estimate}).

By Steps 1 and 2 we know that $E=\xi+\nabla p+ W \in H(\curl,\Om_b)$ is a solution
to the problem (\ref{BVP})-(\ref{BV2}) with the estimate
\be\label{estimate1}
||E||_{H(\curl,\Om_b)}\leq ||\xi||_{H(\curl,\Om_b)}+||\nabla p||_{H(\curl,\Om_b)}
+ ||W||_{H(\curl,\Om_b)}\leq C||W||_{H(\curl,\Om_b)}.
\en
Recalling that
\ben
||f||_{H^{-1/2}_{\dive}(\G_b)}=\inf\{||W||_{H(\curl,\Om_b)}\,|\,
\nu\times W=0\;\mbox{on}\;\G_0\;\mbox{and}\;\nu\times W=f\;\mbox{on}\;\G_b\},
\enn
it follows from (\ref{estimate1}) that
$\ds||E||_{H(\curl,\Om_b)}\leq C||f||_{H^{-1/2}_{\dive}(\G_b)}$.
\end{proof}

For $\ds f\in H^{-1/2}_{\dive}(\G_b)$ define the operator $T$ by
\ben
T(f)=\nu\times(\curl E\times\nu)\qquad{\rm on}\;\;\G_b,
\enn
where $E$ solves the quasi-periodic boundary value problem (\ref{BVP})-(\ref{BV2}).
%with $f\in H^{-1/2}_{\dive}(\G_b)$.
By Lamma \ref{wellpose1}, the operator $T$ is well-defined. Note that $T(f)$ belongs
to the dual space $\ds(H^{-1/2}_{\dive}(\G_b))'=H^{-1/2}_{\curl}(\G_b)$ of
$\ds H^{-1/2}_{\dive}(\G_b)$ with the duality defined by
\ben
<T(f),g>=\int_{\Om_b}[\curl E\cdot\curl\overline{V}-k^2qE\cdot\overline{V}]dx
\enn
for $g\in H^{-1/2}_{\dive}(\G_b),$ where $V\in H(\curl,\Om_b)$ satisfies that
$\nu\times V=g$ on $\G_b$ and $\nu\times V=0$ on $\G_0$.
The operator $T$ can be considered as a Dirichlet-to-Neumann map associated with
the problem (\ref{BVP})-(\ref{BV2}) and depending on the index $q(x)$.
Under the assumptions (A1)-(A3), the above definition of $T(f)$ is independent of the
choice of $V$ and therefore $T: H^{-1/2}_{\dive}(\G_b)\rightarrow
(H^{-1/2}_{\dive}(\G_b))'=H^{-1/2}_{\curl}(\G_b)$ is well-defined.
Moreover, it follows from the above equality and Lemma \ref{wellpose1} that
\ben
||T(f)||_{H^{-1/2}_{\curl}(\G_b)}\leq C||E||_{H(\curl,\Om_b)}
\leq C||f||_{H^{-1/2}_{\dive}(\G_b)}.
\enn
This implies that $T$ is continuous from $H^{-1/2}_{\dive}(\G_b)$ to $H^{-1/2}_{\curl}(\G_b)$.

\section{Solvability of the scattering problem}\label{sec4}
\setcounter{equation}{0}

In this section we will establish the solvability of the scattering problem
(\ref{equation1})-(\ref{RE}), employing the variational method.
To this end, we propose a variational formulation of the scattering problem in
a truncated domain by introducing a transparent boundary condition on $\G_b$.
The existence and uniqueness of solutions to the problem will then be proved using
the Hodge decomposition together with the Fredholm alternative.
%Since this process is almost as the same as the former Lemma,
%we only sketch the proof and refer to \cite{Habib} and \cite{Songlei} for more details.

\subsection{Transparent boundary condition and variational formulation}

Let $x^\prime=(x_1,x_2,b)\in\G_b$ for $b>0$.
For $\widetilde{E}\in H_t^{-\frac12}(\dive,\G_b)$ with
$\widetilde{E}(x^\prime)=\sum_{n\in\Z^2}\widetilde{E}_n\exp(i\al_n\cdot x^\prime)$,
define $\mathcal{R}:\,H_t^{-\frac12}(\dive,\G_b)\rightarrow H_t^{-\frac12}(\curl,\G_b)$
by
\be\label{DtoN}
(\mathcal{R}\widetilde{E})(x^\prime)=(e_3\times\curl E)\times e_3\quad{\rm on}\;\G_b,
\en
where $E$ satisfying the Rayleigh expansion condition (\ref{RE}) is the unique
quasi-periodic solution to the problem
\ben
\curl\curl E-k^2E=0\quad{\rm for}\;x_3>b,\qquad
\nu\times E=\widetilde{E}(x^\prime)\quad{\rm on}\;\G_b.
%E(x)&=&\sum_{n\in\Z^2}E_n\exp(i(\al_n\cdot x^\prime+\beta_nx_3)),\quad .
\enn
The map $\mathcal{R}$ is well-defined and can be used to replace the radiation condition
(\ref{RE}) on $\G_b$. Then the direct scattering problem (\ref{equation1})-(\ref{RE}) can be
transformed into the following boundary value problem in the truncated domain $\Om_b$:
\be\label{equation01}
\curl\curl E-k^2qE&=&0\qquad{\rm in}\;\Om_b,\\ \label{equation02}
\nu\times E&=&0\qquad {\rm on}\;\G_0,\\ \label{equation03}
(\curl E)_T-\mathcal{R}(e_3\times E)&=&(\curl E^i)_T-\mathcal{R}(e_3\times E^i)\quad{\rm on}\;\G_b,
\en
where, for any vector function $U$, $U_T=(\nu\times U)\times\nu$ denotes its
tangential component on a surface. The variational formulation for the problem
(\ref{equation01})-(\ref{equation03}) can be given as follows:
find $E\in X:=\{E\in H(\curl,\Om_b)\;|\;\nu\times E=0\;{\rm on}\;\G_0\}$ such that
\be\label{B}
B(E,\varphi)&:=&\int_{\Om_b}\left[\curl E\cdot\curl\overline{\varphi}
   -k^2qE\cdot\overline{\varphi}\right]dx
   -\int_{\G_b}\mathcal{R}(e_3\times E)\cdot(e_3\times\overline{\varphi})ds\\ \nonumber
   &=&\int_{\G_b}\left[(\curl E^i)_T
     -\mathcal{R}(e_3\times E^i)\right]\cdot(e_3\times\overline{\varphi})ds
\en
for all $\varphi\in X$.
%To study the solvability of (DP),
We have the following properties of $\mathcal{R}$:
\begin{description}
\item 1) $\mathcal{R}\,:\,H_t^{-\frac12}(\dive,\G_b)\rightarrow H_t^{-\frac12}(\curl,\G_b)$
is continuous and can be explicitly represented as
\be\label{p1}
(\mathcal{R}\widetilde{E})(x^\prime)=-\sum_{n\in\Z^2}\frac{1}{i\beta_n}
\left[k^2\widetilde{E}_n-(\al_n\cdot\widetilde{E}_n)\al_n\right]\exp(i\al_n\cdot x^\prime).
\en
\item 2) Let $P=\{n=(n_1,n_2)\in\Z^2\,|\,\beta_n\;\mbox{is a real number}\}$. Then
\be\label{p2}
\Rt<\mathcal{R}\widetilde{E},\widetilde{E}>&=& 4\pi^2\sum_{n\in\Z^2\backslash P}
     \frac{1}{|\beta_n|}\left[k^2|\widetilde{E}_n|^2
        -|\al_n\cdot\widetilde{E}_n|^2\right],\\ \label{p3.1}
-\Rt<\mathcal{R}\widetilde{E},\widetilde{E}>&\geq&
      C_1||\dive\widetilde{E}||^2_{H_t^{-1/2}(\G_b)}
      -C_2||\widetilde{E}||^2_{H_t^{-1/2}(\G_b)},
\en
where $C_1$ and $C_2$ are positive constants and $<\cdot,\cdot>$ denotes
the inner product of $L^2_t(\G_b)$.
\item 3)
\be\label{p3}
\I<\mathcal{R}\widetilde{E},\widetilde{E}>=4\pi^2\sum_{n\in P}\frac{1}{\beta_n}
     \left[k^2|\widetilde{E}_n|^2-|\al_n\cdot\widetilde{E}_n|^2\right]\geq 0.
\en
\end{description}

The representation (\ref{p1}) of $\mathcal{R}$ can be computed directly from its
definition (\ref{DtoN}) (see \cite{AT}) and the properties (\ref{p2})-(\ref{p3})
can be easily obtained using this representation.
Furthermore, there exists a $C>0$ such that for every $\eta>0$ and $E\in H(\curl,\Om_b)$
we have (see \cite{AT})
\be\label{p4}
||\nu\times E||_{H^{-1/2}_t(\G_b)}\leq C
\left[\eta||\curl E||_{L^2(\Om_b)}+(1+1/\eta)||E||_{L^2(\Om_b)}\right].
\en
Let
\ben
S&=&\{p\in H^1(\Om_a)\;|\;p=0\;{\rm on}\;\G_0\}\\
X_0&=&\{E\in X\;|\;B(E,\nabla p)=0\;\forall p\in S\}
%\int_{\Om_b}k^2qE\cdot\nabla\overline{ p}dx
%+\int_{\G_b}\mathcal{R}(e_3\times E)\cdot(e_3\times\nabla\overline{ p})ds=0,\ \forall p\in S.\}
\enn
Then in a completely similar manner as in the proof of Lemma \ref{wellpose},
we can establish the Hodge decomposition $X=X_0\oplus\bigtriangledown S$.

\subsection{Solvability of the direct scattering problem}

\begin{lemma}\label{Strong-ellipticity}
The bilinear form $B(\cdot,\cdot)$ defined by $(\ref{B})$ is strongly
elliptic on $X_0$, that is, for all $w_0\in X_0$,
\ben
\Rt B(w_0,w_0)\geq C||w_0||_{X}-\rho(w_0,w_0)
\enn
for some constant $C>0$ and a compact bilinear form $\rho(\cdot,\cdot)$.
\end{lemma}

\begin{proof}
Let $M$ be a positive constant to be determined later and let
\ben
b_1(w_0,\varphi_0)&=&\int_{\Om_b}\left[\curl w_0\cdot\curl\overline{\varphi_0}
  +(M-k^2q)w_0\cdot\overline{\varphi_0}\right]dx
  -\int_{\G_b}\Rc(e_3\times w_0)\cdot(e_3\times\overline{\varphi}_0)ds\\
b_2(w_0,\varphi_0)&=&-M\int_{\Om_b}w_0\cdot\overline{\varphi_0}dx.
\enn
Then $B(w_0,\varphi_0)=b_1(w_0,\varphi_0)+b_2(w_0,\varphi_0)$ for $w_0,\,\varphi_0\in X_0.$
By the properties of $\Rc$ it follows that
\ben
&&-\Rt<\Rc(e_3\times w_0),e_3\times\overline{w}_0>\\
&&\qquad\ge C_1||\dive(e_3\times w_0)||^2_{H_t^{-1/2}(\G_b)}
            -C_2||e_3\times w_0 ||^2_{H_t^{-1/2}(\G_b)}\\
&&\qquad\ge C_1||\dive(e_3\times w_0)||^2_{H_t^{-1/2}(\G_b)}
         -C_3\eta^2||\curl w_0||^2_{L^2(\Om_b)}
         -C_3(1+\frac1{\eta})^2||w_0||^2_{L^2(\Om_b)},
\enn
where $C_1,\,C_2$ and $C_3$ are three positive constants and $\eta>0$ is arbitrary.
Thus we have
\ben
\Rt b_1(w_0,w_0)&\ge&\|\curl w_0\|^2_{L^2(\Om_b)}+(M-k^2q_\infty)\|w_0\|^2_{L^2(\Om_b)}\\
      &&-C_3\eta^2||\curl w_0||^2_{L^2(\Om_b)}-C_3(1+1/\eta)^2||w_0||^2_{L^2(\Om_b)}\\
   &=&(1-C_3\eta^2)||\curl w_0||^2_{L^2(\Om_b)}+[M-k^2q_\infty-C_3(1+1/\eta)^2]||w_0||^2_{L^2(\Om_b)},
\enn
where $q_\infty=\max_{x\in\R^3_+}|q(x)|<\infty.$
Choose $\eta$ sufficiently small and $M$ sufficiently large so that
\be\label{3.2}
\Rt b_1(w_0,w_0)\geq C_0(||\curl w_0||^2_{L^2(\Om_b)}+||w_0||^2_{L^2(\Om_b)})
\en
for some constant $C_0>0$. This, together with the fact that $X_0$ is compactly
imbedded in $(L^2(\Om_b))^3$, yields the desired result.
\end{proof}

\begin{theorem}\label{wellpose}
Assume that the conditions $(A1)-(A3)$ are satisfied. Then the problem
$(\ref{equation1})-(\ref{RE})$ has a unique solution $E\in H_{loc}(\curl,\R^3_+)$ such that
\ben\label{estimate2}
||E||_{H_{loc}(\curl,\R^3_+)}:=\max_{a>b}||E||_{H(\curl,\Om_a)}\leq C||E^i||_{H(\curl,\Om_b)},
\enn
where $C$ is a positive constant depending on the domain and $q$.
\end{theorem}

\begin{proof}
It follows from Lemma \ref{Strong-ellipticity} and the proof of Lemma \ref{wellpose1}
that there exists a unique solution $E\in H(\curl,\Om_b)$ satisfying that
$||E||_{H(\curl,\Om_b)}\leq C||E^i||_{H(\curl,\Om_b)}$.
It remains to extend $E(x)$ to be a function in $H_{loc}(\curl,\mathbb{R}^3_+)$.
Suppose
$e_3\times(E-E^i)|_{\G_b}=\sum_{n\in\mathbb{N}\times\mathbb{N}}
A_n e^{i\al_n\cdot x}\in H^{-1/2}(\dive,\G_b)$.
Let
\ben
E^s(x)=\sum_{n\in\mathbb{N}\times\mathbb{N}}(A_n\times e_3+B_ne_3)
e^{i\al_n\cdot x+i\beta_n(x_3-b)},\qquad x_3>b
\enn
and let $E^s$ satisfy that $\dive E^s(x)=0$ for $x_3>b$. Then we have
$B_n=\frac{1}{\beta_n}(e_3\times A_n)\cdot\al_n$. Thus
\ben
E^s(x)=\sum_{n\in \mathbb{N}\times\mathbb{N}}\left[A_n\times e_3
+\frac{1}{\beta_n}(e_3\times A_n)\cdot\al_ne_3\right]
e^{i\al_n\cdot x+i\beta_n(x_3-b)},\qquad x_3>b.
\enn
Define $E(x)=E^i(x)+E^s(x)$ for $x_3>b$. Then it is easy to prove that
$E\in H(\curl,\Om_a\ba\Om_b)$ with
$||E||_{H(\curl,\Om_a\ba\Om_b)}\leq C||E^i||_{H(\curl,\Om_b)}$
for any $a>b$, so $E\in H(\curl,\Om_a)$ for any $a>b$, that is,
$E\in H_{loc}(\curl,\mathbb{R}^3_+)$ with the required estimate (\ref{estimate2}).
The proof is thus completed.
%Since $q(x)\equiv 1$ in a neighborhood of $\G_b$, the application of standard
%elliptic theory leads to that $E(x)$ is smooth up to  $\G_b$. This implies that
%the previously  defined $E(x)$ belongs to $H(curl,\Om_a)$ for any $a>b$.
\end{proof}

\section{The inverse problem}\label{inverse problems}
\setcounter{equation}{0}

Let $a>b$ and assume that there are two refractive index functions $q_i\;(i=1,2)$
satisfying the assumptions (A1)-(A3). For $g\in L_t^2(\G_a)$ let the incident waves
be of the form:
\be\label{incident waves}
E^i(x,g)=\curl_x\curl_x\int_{\G_a}G(x,y)g(y)ds(y),\qquad x<a.
\en
Write the scattered electric field and the total electric field as
$E^s_i(x,g)$ and $E_i(x,g)$, respectively, indicating their dependance on $g$ and
the refractive index function $q_i$ $(i=1,2).$

For the refractive index $q_i$ denote by $T_i$ the corresponding Dirichlet-to-Neumann
map associated with the problem (\ref{BVP})-(\ref{BV2}) with $q$ replaced by $q_i$ $(i=1,2)$,
as defined at the end of Section \ref{sec3}.
%
%In view of the D-to-N map $T$ which was introduced at the end of Section \ref{sec3}
%and associated with the QPBVP (\ref{BVP})-(\ref{BV2}), we denote by $T_i(i=1,2)$ the
%corresponding Dirichlet to Neumann map depending on the index $q_i(x)$.

\begin{lemma}\label{T}
If $T_1(f)=T_2(f)$ for all $f\in H_t^{-1/2}(\dive,\G_b)$, then
\ben
\int_{\Om_b} E_1(x)\cdot\overline{E}_2(x)\left[q_1(x)-q_2(x)\right]dx=0,
\enn
where $E_1,\,E_2\in H(\curl,\Om_b)$ solve the problem $(\ref{BVP})-(\ref{BV2})$ with $q$
replaced by $q_1$ and $\overline{q}_2$, respectively.
\end{lemma}

\begin{proof}
Let $E_1$ and $F_2\in H(\curl,\Om_b)$ be the solution of the problems
\ben
\curl\curl E_1-k^2q_1E_1=0\quad{\rm in}\;\Om_b,\qquad\nu\times E_1=0\quad{\rm on}\;\G_0
\enn
and
\ben
\curl\curl F_2-k^2q_2F_2=0\quad{\rm in}\;\Om_b,\quad\nu\times F_2=0\quad{\rm on}\;\G_0,
\quad\nu\times F_2=\nu\times E_1\quad{\rm on}\;\G_b,
\enn
respectively. Let $E=F_2-E_1$. Then it is easy to see that
\ben
\curl\curl E-k^2q_2E&=&k^2(q_2-q_1)E_1\quad{\rm in}\;\Om_b,\\
\nu\times E&=&0\quad{\rm on}\;\G_0\cup\G_b,\\
\nu\times\curl E&=&0\quad{\rm on}\;\G_b,
\enn
where the last quality is obtained from the assumption $T_1=T_2$.
Thus, it follows from the Green vector formula that
\ben
\int_{\Om_b}(q_2-q_1)E_1\cdot\overline{E}_2dx
&=&\frac{1}{k^2}\int_{\Om_b}(\curl\curl E-k^2q_2E)\cdot\overline{E}_2dx\\
&=&\frac{1}{k^2}\int_{\Om_b}(\curl E\cdot\curl\overline{E}_2-k^2q_2E\cdot\overline{E}_2)dx\\
&=&\frac{1}{k^2}\int_{\Om_b}(E\cdot\curl\curl\overline{E}_2-k^2q_2E\cdot\overline{E}_2)dx\\
&=&\frac{1}{k^2}\int_{\Om_b}(E\cdot k^2q_2\overline{E}_2-k^2q_2E\cdot\overline{E}_2)dx=0.
\enn
The proof is thus completed.
\end{proof}

For $g\in L_t^2(\G_a)$ appearing in the incident waves (\ref{incident waves}), we define an
operator $F:L_t^2(\G_a)\rightarrow H_t^{-1/2}(\dive,\G_b)$ by
\ben
F(g)=e_3\times E(x,g)\quad{\rm on}\quad\G_b,
\enn
where $E(x,g)$ solves the problem (\ref{equation1})-(\ref{equation4}) with the incident wave
$E^i(x,g)$. The operator $F$ can be considered as an input-output operator mapping the sum of
the electric dipoles to the tangential component of the corresponding total field on $\G_b$.
Moreover, for all $g\in L_t^2(\G_a)$, the operator $F$ has a dense range in $H_t^{-1/2}(\dive,\G_b)$,
as stated in the following lemma.

\begin{lemma}\label{dense}
The operator $F$ has a dense range in $H_t^{-1/2}(\dive,\G_b).$
\end{lemma}

\begin{proof}
We only need to prove that $F^*:H_t^{-1/2}(\curl,\G_b)\rightarrow L_t^2(\G_a)$ is injective.
First, we show that for any $f\in H_t^{-1/2}(\curl,\G_b)$, $F^*(f)$ is given by
\be\label{con}
F^*(f)=\left[\curl_y\curl_y\int_{\G_b}\overline{G(x,y)}
\curl(\overline{V^+(x)-W(x)})\times e_3ds(x)\right]_T,
\en
where the superscripts $+$ and $-$ indicate the limit obtained
from $\R_3\ba\Om_b$ and $\Om_b$, respectively,
%$U_T=(\nu\times U)\times\nu$ denotes its tangential component on a
%surface for any vector function $U$,
and for any $a>b$ the function $V\in H(\curl,\Om_b)\cap H(\curl,\Om_a\ba\Om_b)$
solves the problem
\be
\curl\curl V-k^2V=0 &&\mbox{for}\quad x_3>b,\\ \label{equation for V}
\curl\curl V-k^2qV=0  &&\mbox{in}\quad\Om_b,\\ \label{V=0}
\nu\times V=0       &&\mbox{on}\quad\G_0,\\ \label{transmission1}
\nu\times V^+-\nu\times V^-=0\ &&\mbox{on} \quad\G_b,\\ \label{transmission2}
\left[\curl V^+-\curl V^-\right]_T=\overline{f}&&\mbox{on} \quad\G_b
\en
and satisfies the Rayleigh expansion condition (\ref{RE}) with $\al$ replaced
by $-\al$ for $x_3>b$, that is,
\be\label{REV}
V(x)=\sum_{n\in\Z^2}V_ne^{i(\al_n^{\prime}\cdot x+\bt_n^{\prime}x_3)},\qquad x_3\geq b
\en
with $\al_n^{\prime}=(-\al_1+n_1,-\al_2+n_2,0)\in\R^3$, $V_n\in\C^3$ and
$$
\beta_n^{\prime}=\left\{\begin{array}{lll}
             (k^2-|\al_n^{\prime}|^2)^{\frac{1}{2}}\qquad\rm{if}\ |\al_n^{\prime}|< k,\\
              i(|\al_n^{\prime}|^2-k^2)^{\frac{1}{2}}\qquad\rm{if}\ |\al_n^{\prime}|> k.
             \end{array}\right.
$$
In addition, the function $W$ is given by
\be\label{W}
W(x)=\sum_{n\in\Z^2}V_ne^{i((\al_n^{\prime}\cdot x+\bt_n^{\prime}(2b-x_3))},
\qquad x_3\leq b.
\en
In fact, for any $f\in H_t^{-1/2}(\curl,\G_b)$ and $g\in H_t^{-1/2}(\dive,\G_b)$ we have
\ben\no
&&<Fg,f>_{H_t^{-1/2}(\dive,\G_b)\times H_t^{-1/2}(\curl,\G_b)}\\ \no
&&\quad=\int_{\G_b}\nu\times E(\cdot,g)\cdot\overline{f}ds\\ \no
&&\quad=\int_{\G_b}\nu\times E(\cdot,g)\cdot[\curl V^+-\curl V^-]ds\\
&&\quad=\int_{\G_b}[(\nu\times E\cdot\curl V^+-\nu\times V^+\cdot\curl E)
        -(\nu\times E\cdot\curl V^--\nu\times V^-\cdot\curl E)]ds,
\enn
where the transmission conditions (\ref{transmission1}) and (\ref{transmission2})
have been used. It follows from the Maxwell equations (\ref{equation for V}) and
(\ref{equation2}) and the boundary conditions (\ref{V=0}) and (\ref{equation3}) that
\be\label{F2}
\int_{\G_b}[\nu\times E\cdot\curl V^--\nu\times V^-\cdot\curl E]ds=0.
\en
On the other hand, from the Rayleigh expansion conditions (\ref{RE}) and (\ref{REV})
it is derived that
\be\no
&&\int_{\G_b}[\nu\times E\cdot\curl V^+-\nu\times V^+\cdot\curl E]ds\\ \no
&&\quad=\int_{\G_b}[(\nu\times E^i\cdot\curl V^+-\nu\times V^+\cdot\curl E^i)
   +(\nu\times E^s\cdot\curl V^+-\nu\times V^+\cdot\curl E^s)]ds\\ \label{F3}
&&\quad=\int_{\G_b}[\nu\times E^i\cdot\curl V^+-\nu\times V^+\cdot\curl E^i]ds.
\en
Similarly, from the definition of $E^i$ and the Rayleigh expansion condition (\ref{W})
it follows that
\be\label{F4}
\ds\int_{\G_b}[\nu\times E^i\cdot\curl W-\nu\times W\cdot\curl E^i]ds=0.
\en
The equations (\ref{F2})-(\ref{F4}) together with the fact that $V=W$ on $\G_b$ yield
\ben
<Fg,f>&=&\int_{\G_b}[\nu\times E^i\cdot\curl V^+-\nu\times V^+\cdot\curl E^i]ds\\
&=&\int_{\G_b}[\nu\times E^i\cdot\curl V^+-\nu\times W\cdot\curl E^i]ds\\
&=&\int_{\G_b}[\nu\times E^i\cdot\curl V^+-\nu\times E^i\cdot\curl W]ds\\
&=&\int_{\G_b}\nu\times E^i\cdot(\curl V^+-\curl W)ds.
\enn
Substituting the expression (\ref{incident waves}) of $E^i$ into the above equation
and exchanging the order of integration we get
\ben
<Fg,f>=\int_{\G_a}g(y)\cdot\curl_y\curl_y\left[\int_{\G_b}G(x,y)\curl[V^+(x)-W(x)]
     \times e_3ds(x)\right]ds(y),
\enn
which implies (\ref{con}).

We now prove that $F^*$ is injective. Suppose $F^*(f)=0$ for some
$f\in H_t^{-1/2}(\curl,\G_b)$. Define $U$ by
\ben
U(y):=\curl_y\curl_y\left[\int_{\G_b}\overline{G(x,y)}h(x)ds(x)\right],
\qquad y\in\mathbb{R}^3\ba\G_b,
  \enn
where $h=\curl(\overline{V^+-W})\times e_3$. Then $e_3\times U(y)=0$ on $\G_a$.
It is clear that $U(y)$ is a $-\al$-quasi-periodic function satisfying the Rayleigh expansion
condition (\ref{RE}) when $y_3>a$. By the uniqueness of solutions to the exterior Dirichlet
problem (see \cite{AH}) we have $U(y)=0$ when $y_3>a$, which together with
the unique continuation principle (\cite{CK}) implies that $U(y)=0$ when $y_3>b.$
Now from the jump relation $e_3\times U^+(y)-e_3\times U^-(y)=0$ on $\G_b$
and again the uniqueness of solutions for the exterior Dirichlet problem for $y_3<b$
we get that $U(y)=0$ when $y_3<b$. Thus, $h(y)=e_3\times\curl[U^+(y)-U^-(y)]=0$ on $\G_b,$
which, together with (\ref{REV}) and (\ref{W}), implies that
\be\label{wv}
e_3\times V^+=e_3\times W,\quad e_3\times\curl V^+=e_3\times\curl W\quad{\rm on}\;\G_b.
\en
Since $V$ and $W$ satisfy the Maxwell equation $\curl\curl E-k^2E=0$ in the regions
$x_3>b$ and $x_3<b$, respectively, then it follows easily from the transmission
condition (\ref{wv}) and the Rayleigh expansion conditions (\ref{REV}) and (\ref{W})
that $V=0$ for $x_3>b$ and $W=0$ for $x_3<b.$ Thus, by (\ref{transmission1})
we have $\nu\times V^-=0$ on $\G_b,$ so $V\in H(\curl,\Om_b)$ satisfies
the problem (\ref{BVP})-(\ref{BV2}) with $f=0$. By Lemma \ref{wellpose1}
we have $V=0$ in $\Om_b$.
%$V$ satisfies upward $-\al$ Rayleigh expansion and $W$ satisfies down-ward
%$-\al$ Rayleigh expansion, we have $V=W=0$ in $\mathbb{R}^3$ which yield
Thus, $f=[\curl\overline{V}^+-\curl\overline{V}^-]_T=0$,
which completes the proof of Lemma \ref{dense}.
\end{proof}

Combining Lemmas \ref{T} and \ref{dense}, we have the following orthogonality
relation for two different functions $q_i\;(i=1,2)$.

\begin{lemma}\label{orthogonality}
Let the incident waves $E^i(x,g)$ be defined by $(\ref{incident waves}).$
If
\be\label{orth}
e_3\times E_1(x,g)=e_3\times E_2(x,g)\qquad{\rm on}\;\;\G_a
\en
for all $g\in L_t^2(\G_a)$ and some $a>b$, then the following orthogonality relation holds:
\ben
\int_{\Om_b}E_1(x)\cdot\overline{E}_2(x)(q_1(x)-q_2(x))dx=0,
\enn
where $E_1,\;E_2\in H(\curl,\Om_b)$ solve the problem $(\ref{BVP})-(\ref{BV2})$
with $q$ replaced by $q_1$ and $\overline{q}_2$, respectively.
\end{lemma}

\begin{proof}
From the equation (\ref{orth}), the uniqueness of solutions for the
exterior Dirichlet problem and the unique continuation principle it follows that
$E_1(x,g)= E_2(x,g)$ for all $x_3>b$. This implies that
\ben
e_3\times\curl E_1^+(x,g)=e_3\times\curl E_2^+(x,g)\qquad{\rm on}\;\;\G_b.
\enn
Since $[e_3\times\curl E_j^+(x;g)]|_{\G_b}=0$ for $j=1,2$, then we have
\ben
e_3\times\curl E_1^-(x,g)=e_3\times\curl E_2^-(x,g)\qquad{\rm on}\;\;\G_b.
\enn
%Combining the above two equalities, in view of (\ref{BVP}) and definitions of $T_i$ we get
By the above two equalities and the definition of $T_i$ we have
\ben
T_1(e_3\times E_1(x,g))=T_2(e_3\times E_2(x,g))
\enn
for all $g\in L_t^2(\G_a).$ The continuity of $T_j\;(j=1,2)$ and Lemma \ref{dense} lead to
\ben
T_1(f)=T_2(f)\qquad\forall f\in H_t^{-1/2}(\dive,\G_b).
\enn
This together with Lemma \ref{T} gives the desired result.
\end{proof}

We are now ready to prove our main result for the inverse scattering problem.

\begin{theorem}\label{thm-inverse}
Let $q_j\;(j=1,2)$ satisfy the assumptions $(A1)-(A3)$ and let $q_j$ depend on
only one direction $x_1$ or $x_2$ with $j=1,2$. If
\ben
e_3\times E_1(x,g)=e_3\times E_2(x,g)\quad{\rm on}\;\;\G_a
\enn
for all $g\in L_t^2(\G_a)$ with some $a>b$, where $E_j(x,g)$ solves the problem
$(\ref{equation1})-(\ref{equation4})$ with $q=q_j$ $(j=1,2)$ corresponding to
the incident wave $E^i(x,g)$ given by $(\ref{incident waves})$, then $q_1=q_2$.
\end{theorem}

\begin{proof}
By Lemma $\ref{orthogonality}$ we have the orthogonality relation:
\be\label{orth2}
\int_{\Om_b}E_1(x)\cdot\overline{E}_2(x)\left[q_1(x)-q_2(x)\right]dx=0,
\en
where $E_1,\;E_2\in H(\curl,\Om_b)$ solve the problem (\ref{BVP})-(\ref{BV2}) with $q$
replaced by $q_1$ and $\overline{q}_2$, respectively.

We now look for solutions to the problem (\ref{BVP})-(\ref{BV2}) in the following form:
\ben
E(x)=(0,0,E_3(x_1,x_2))=(0,0,v(x_1)u(x_2))
\enn
with the scalar functions $v$ and $u$ satisfying the following quasi-periodic conditions:
\ben
v(x_1)e^{2i\al_1\pi}=v(x_1+2\pi),\qquad u(x_2)e^{2i\al_2\pi}=v(x_2+2\pi).
\enn
It is clear that such a function $E$ is $\al$-quasi-periodic and satisfies the
boundary condition (\ref{BV1}).
Without loss of generality, we may assume that $q_j(x)=q_j(x_1)$,
that is, $q_j$ depends only the $x_1$-direction with $j=1,2$.
Substituting such $E$ into the Maxwell equation (\ref{BVP}) and noting that
$\curl\curl=-\triangle+\nabla(\nabla\cdot)$, we find that
\ben
v''(x_1)u(x_2)+v(x_1)u''(x_2)+k^2q(x_1)v(x_1)u(x_2)=0,\qquad x_1,\,x_2\in(0,2\pi),
\enn
which implies that
\ben
\frac{v''(x_1)}{v(x_1)}+k^2q(x_1)v(x_1)=\frac{u''(x_2)}{u(x_2)}=\lambda
\enn
for some constant $\lambda$, where $x_1,x_2\in(0,2\pi).$
Following the idea of Kirsch \cite{Kirsch95}, we construct a special kind of solutions
$v$ by considering the following quasi-periodic Sturm-Liouville eigenvalue problem:
\ben
(\textrm{I}):\quad\left\{\begin{array}{lll}
                   v''(x_1)+k^2q(x_1)v(x_1)=\lambda v(x_1),\quad x_1\in(0,2\pi)\\
                   v(x_1)e^{2i\al_1\pi}=v(x_1+2\pi),\\
                   v'(x_1)e^{2i\al_1\pi}=v'(x_1+2\pi).
                 \end{array}\right.
\enn
The eigenvalues $\lambda_n$ and the corresponding eigenfunctions $v_n$, normalized to $v_n(0)=1$,
have the following asymptotic behaviors as $n\rightarrow\infty$ (see \cite{W}):
\ben
\lambda^{\pm}_n&=&\left(n\pm\frac{\al_1}{2\pi}\right)^2
   -\frac{k^2}{2\pi}\int_0^{2\pi}q(s)ds+\mathcal{O}\left(\frac{1}{n}\right),\\
v_n^{\pm}(x_1)&=&\exp\left[{i(\pm n+\frac{\al_1}{2\pi})x_1}\right]
+\mathcal{O}\left(\frac{1}{n}\right)
\enn
which are uniform in $x_1\in [0,2\pi]$. We also consider the following quasi-periodic
boundary problem for $u$:
\ben
(\textrm{II}):\quad\left\{\begin{array}{lll}
                   u''(x_2)-\lambda_n u(x_2)=0,\quad x_2\in(0,2\pi)\\
                   u(x_2)e^{2i\al_2\pi}=v(x_2+2\pi).
                  \end{array}\right.
\enn
The non-trivial solutions to the problem (\textrm{II}) can be written explicitly as
\ben
u_n(x_2)=c_{n,1}e^{\sqrt{\lambda}_nx_2}+c_{n,1}e^{-\sqrt{\lambda}_nx_2},\quad\lambda_n\neq0,
\enn
where $c_{n,1}$ and $c_{n,2}$ are constants satisfying
\be\label{C}
c_{n,1}=c_{n,2}\left(e^{-2\pi\sqrt{\lambda}_n}-e^{i2\pi\al_2}\right)\Big/
\left(e^{i2\pi\al_2}-e^{2\pi\sqrt{\lambda}_n}\right).
\en

Now, let $E_{3,n}^{\pm}=v_n^{\pm}(x_1)u_n^{\pm}(x_2)$ be the third component of
$E_n^{\pm}=(0,0,E_{3,n}^{\pm})$ corresponding to $q_1(x_1)$ and let
$E_{3,m}^{\pm}=v_m^{\pm}(x_1)u_m^{\pm}(x_2)$ be the third component of $E_n^\pm$
corresponding to $\overline{q_2}(x_1)$. It follows from (\ref{orth2}) that
\be\label{q12}
0=\int_{\Om_b} E_{3,n}(x_1,x_2)\cdot\overline{E}_{3,m}(x_1,x_2)
\left[q_1(x_1)-q_2(x_1)\right]dx=bA_1^{n,m}A_2^{n,m},
\en
where
\ben
A_1^{n,m}:&=&\int_0^{2\pi}[q_1(x_1)-q_2(x_1)]e^{i(n-m)x_1}dx_1
             +\mathcal{O}\left(\frac{1}{n}\right)+\mathcal{O}\left(\frac{1}{m}\right),\\
A_2^{n,m}:&=&\int_0^{2\pi}\left(c_{n,1}e^{\sqrt{\lambda}_nx_2}
             +c_{n,2}e^{-\sqrt{\lambda}_nx_2}\right)
            \left(\overline{c_{m,1}e^{\sqrt{\lambda}_mx_2}
                   +c_{m,2}e^{-\sqrt{\lambda}_mx_2}}\right)dx_2
\enn
and $c_{n,j},\;c_{m,j}$ satisfy (\ref{C}) with $j=1,2.$
For arbitrarily fixed $l\in\mathbb{N}$, letting $m=n-l$ gives
\ben
A_1^{m+l,m}&=&\int_0^{2\pi}[q_1(x_1)-q_2(x_1)]e^{ilx_1}d(x_1)
              +\mathcal{O}\left(\frac{1}{m}\right),\\
A_2^{m+l,m}&=&\int_0^{2\pi}\left(c_{m+l,1}e^{\sqrt{\lambda_{m+l}}x_2}
              +c_{m+l,2}e^{-\sqrt{\lambda_{m+l}}x_2}\right)
             \left(\overline{c_{m,1}e^{\sqrt{\lambda}_mx_2}
             +c_{m,2}e^{-\sqrt{\lambda}_mx_2}}\right)dx_2.
\enn
We can always choose appropriate constants $c_{m,2}$ and $c_{m,1}$
satisfying (\ref{C}) such that $A_2^{m+l,m}\neq 0$ for sufficiently large $m$.
In fact, we may assume that $l$ is a positive number since otherwise we can
take $n=m-l^\prime$ for some positive $l^\prime$ instead of $l$.
Now choose $\ds c_{m,2}=e^{2\pi\sqrt{\lambda_m}}$. Then, by (\ref{C}), $|c_{m,1}|\ge C_1$
for large $m$ with some positive constant $C_1$ independent of $m$ and
$\left|\int_0^{2\pi}c_{m,2}e^{-2\pi\sqrt{\lambda_m}x_2}dx_2\right|$ tends to $+\infty$
as $m\rightarrow\infty$. This implies that $|A_2^{m+l,m}|\rightarrow+\infty$
as $m\rightarrow+\infty$. Letting $m\rightarrow+\infty$ we conclude from (\ref{q12}) and
the above discussion that
$$
\ds\int_0^{2\pi}(q_1(x_1)-q_2(x_1))e^{ilx_1}dx_1=0
$$
for every $l\in\N$, which implies that $q_1=q_2$. The proof is thus completed.
\end{proof}

\section*{Acknowledgements}

This work was supported by the NNSF of China grant No. 10671201.

\end{document}